\documentclass[12pt,leqno]{article}
\usepackage{amsmath,amssymb,amsfonts,amsthm,latexsym,amstext,amscd}
\usepackage[english]{babel}
\usepackage{fancyhdr}
\newtheorem{theorem}{Theorem}
\newtheorem{corollary}[theorem]{Corollary}

\theoremstyle{definition}

\newtheorem*{defi}{Definition}
\newtheorem*{rema}{Remark}

\newcommand{\beq}{\begin{equation}}
\newcommand{\eeq}{\end{equation}}

\def\ve{\varepsilon }
\begin{document}

\title%[A note on the distribution of normalized prime gaps]
{A note on the distribution of normalized prime gaps}

\author%[J. Pintz]
{J\'anos Pintz\thanks{Supported by OTKA Grants NK104183, K100291 and ERC-AdG.~321104.}
}
                      
\date{}

\numberwithin{equation}{section}
%\begin{abstract}\end{abstract}\begin{classification}Primary 11N05, 11N36; Secondary 11N35.\end{classification}

%\begin{keywords}Hardy--Littlewood prime $k$-tuple conjecture, prime numbers, sieves, gaps between primes, twin primes.
%\end{keywords}

\maketitle

\section{Introduction}
\label{sec:1}

The Prime Number Theorem implies that the average value of
\beq
\label{eq:1.1}
d_n = p_{n + 1} - p_n
\eeq
is $(1 + o(1)) \log p_n$ if $n \in [N, 2N]$, for example, where $\mathbb P = \{p_i \}_{i = 1}^\infty$ is the set of primes.
This motivates the investigation of the sequence $\{d_n / \log p_n\}_{n = 1}^\infty$ or $\{d_n / \log n\}_{n = 1}^\infty$ (which is asymptotically equal).
Erd\H{o}s formulated the conjecture that the set of its limit points
\beq
\label{eq:1.2}
J = \left\{\frac{d_n}{\log n} \right\}' = [0, \infty].
\eeq
He writes in \cite{Erd1955}:
``It seems certain that $d_n / \log n$ is everywhere dense in $(0, \infty)$'' (after mentioning the conjecture $\liminf\limits_{n \to \infty} d_n / \log n = 0$).
The fact that $\infty \in J$ was proved already in 1931 by Westzynthius \cite{Wes1931}.

In 2005 Goldston, Y{\i}ld{\i}r{\i}m and the author \cite{GPY2006}, \cite{GPY2009} showed
$0 \in J$ which is the hitherto only concrete known element of~$J$.
On the other hand already 60 years ago Ricci \cite{Ric1954} and Erd\H{o}s \cite{Erd1955} proved (simultaneously and independently) that $J$ has a positive Lebesgue measure.
A partial result towards the full conjecture \eqref{eq:1.2} was shown by the author in \cite{Pin2013arX} according to which there exists an ineffective constant $c$ such that
\beq
\label{eq:1.3}
[0, c] \subset J.
\eeq

In a recent work W. Banks, T. Freiberg and J. Maynard \cite{BFM2014arX} proved that for any sequence of $k = 9$ non\-negative real numbers $\beta_1 \leqslant \beta_2 \leqslant \cdots \leqslant \beta_k$ we have
\beq
\label{eq:1.4}
\{\beta_j - \beta_i :\ 1 \leqslant i < j \leqslant k\} \cap J \neq \emptyset.
\eeq
As a corollary they obtained that if $\lambda$ denotes the Lebesgue measure, then
\beq
\label{eq:1.5}
\lambda \bigl([0, T] \cap J\bigr) \geqslant \bigl(1 + o(1)\bigr)T/8.
\eeq

\section{Generalization and Improvement}
\label{sec:2}

The purpose of this note is to generalize this result for the case when $d_n$ is normalized by a rather general function $f(n) \to \infty$, that is to consider instead of $J$ the more general case of the set of limit points
\beq
\label{eq:2.1}
J_f = \left\{\frac{d_n}{f(n)}\right\}'
\eeq
where we require from $f$ to belong to the class $\mathcal F$ below.

\begin{defi}
A function $f(n) \nearrow \infty$ belongs to $\mathcal F$ if for any $\ve > 0$
\beq
\label{eq:2.2}
(1 - \ve) f(N) \leqslant f(n) \leqslant(1 + \ve) f(N) \ \text{ for } n \in [N, 2N], \ N > N_0,
\eeq
further if
\beq
\label{eq:2.3}
f(n) \ll \log n \log_2 n \log_4 n \Big/ (\log_3 n)^2
\eeq
where $\log_\nu n$ denotes the $\nu$-times iterated logarithm.
\end{defi}

The first condition means that $f(n)$ is slowly oscillating, while the second one that it does not grow more quickly than the Erd\H{o}s--Rankin function, which until the recent dramatic new developments by Maynard \cite{May2014arX}, Ford--Green--Konyagin--Tao \cite{FGKT2014arX},
and Ford--Green--Kony\-agin--May\-nard--Tao \cite{FGKMT2014arX} described the largest known gap between consecutive primes.
The improvement means that it is sufficient to work with $k = 5$ values of $\beta_i$ in \eqref{eq:1.4} instead of $k = 9$ values.
As an immediate corollary we obtain a lower bound $(1 + o(1))T/4$ instead of \eqref{eq:1.5} for the Lebesgue measure of the more general set $[0, T] \cap J_f$.

\begin{theorem}
\label{th:1}
If $f \in \mathcal F$, then for any sequence of $k = 5$ nonnegative real numbers $\beta_1 \leqslant \beta_2 \leqslant \dots \leqslant \beta_k$ we have
\beq
\label{eq:2.4}
\bigl\{\beta_j - \beta_i : \ 1 \leqslant i < j \leqslant k\bigr\} \cap J_f \neq \emptyset.
\eeq
\end{theorem}

\begin{corollary}
\label{cor:2}
If $f \in \mathcal F$, then
\beq
\label{eq:2.5}
\lambda\bigl([0, T] \cap J\bigr) \geqslant \bigl(1 + o(1)\bigr)T/4.
\eeq
\end{corollary}

In an earlier work \cite{Pin2013arX} we showed that for any $f \in \mathcal F$ there exists an ineffective constant $c_f$ such that $[0, c_f] \subset J_f$.
We further remark that since $\beta_i$ can be arbitrarily large,
Theorem~\ref{th:1} includes the improvement of the Erd\H{o}s--Rankin function given in \eqref{eq:2.3} proved recently in \cite{May2014arX} and \cite{FGKT2014arX}.
(We note that the proof uses some refinement of the argument of \cite{May2014arX}, so it does not represent an independent new proof.)

In connection with the original Erd\H{o}s conjecture for general $f \in \mathcal F$ we remark that it was proved in \cite{Pin2014arX} that the conjecture is in some sense valid for almost all functions $f \in \mathcal F$.
More precisely it was shown in \cite{Pin2014arX} that if $\{f_n\}_1^\infty \in \mathcal F$ with $\lim\limits_{x \to \infty} f_{n + 1}(x) / f_n(x) = \infty$, then
\beq
\label{eq:2.6}
J_{f_n} = [0, \infty]
\eeq
apart from at most 98 exceptional functions $f_n$.

\section{Proof}
\label{sec:3}

The generalization for the case $f\in \mathcal F$ instead of the single case $f = \log n$ runs completely analogously to the proofs in \cite{Pin2014arX} so we will only describe how to improve $k = 9$ to $k = 5$ in Theorem~\ref{th:1} which leads to the improved Corollary~\ref{cor:2} in the same simple way as described in the Introduction of the work of Banks, Freiberg and Maynard \cite{BFM2014arX}.

The result will follow from the following improvement of Theorem 4.3 of their work.
Let $\mathcal Z$ be given by (4.8) of \cite{BFM2014arX}.

\begin{theorem}
\label{th:3}
Let $m, k$ and $\ve  = \ve (k)$ be fixed.
If $k$ is a sufficiently large multiple of $4m + 1$ and $\ve $ is sufficiently small, there is some $N(m, k, \ve )$ such that the following holds for $N \geqslant N(m, k, \ve )$ with
\beq
\label{eq:3.1}
w = \ve  \log N, \ \ \ W = \prod_{p \leqslant w, p \nmid \mathcal Z} p.
\eeq

Let $\mathcal H = \left\{h_1, \dots, h_k\right\}$ be an admissible $k$-tuple
(that is it does not cover all residue classes $\text{\rm mod }p$ for any prime~$p$) such that
\beq
\label{eq:3.2}
0 \leqslant h_1 < \dots < h_k \leqslant N
\eeq
and
\beq
\label{eq:3.3}
p \bigm|  \prod_{1 \leqslant i < j \leqslant k} (h_j - h_i) \Longrightarrow p \leqslant w.
\eeq
Let $\mathcal H = \mathcal H_1 \cup \dots \cup \mathcal H_{4n + 1}$ be a partition of $\mathcal H$ into $4m + 1$ sets of equal size and let $b$ be an integer with
\beq
\label{eq:3.4}
\biggl(\prod_{i = 1}^k (b + h_i), W\biggr) = 1.
\eeq
Then there is some $n \in (N, 2N]$ with $n \equiv b(\text{\rm mod }W)$ and some set of distinct indices $\{i_1, \dots, i_{m + 1}\} \subseteq \{1, \dots, 4m + 1\}$ such that
\beq
\label{eq:3.5}
\bigl|\mathcal H_i(n) \cap \mathbb P \bigr| \geqslant 1 \ \ \text{ for all } i \in \{i_1, \dots, i_{m + 1}\}.
\eeq
\end{theorem}

\begin{rema}
The original analogous statement (4.20) of \cite{BFM2014arX} should have been stated with $\geqslant 1$ instead of $=1$ (oral communication of James Maynard).
This form is enough to imply their Corollary 1.2 or our Corollary~\ref{cor:2}.
\end{rema}

The needed change in the Deduction of Theorem 4.3 is the following.
First, using $4m + 1 \mid k$ we write
\beq
\label{eq:3.6}
\mathcal H = \mathcal H_1 \cup \dots \cup \mathcal H_{4m + 1}
\eeq
as a partition of $\mathcal H$ into $4m + 1$ sets each of size $k/(4m + 1)$.
Instead of the quantity $S$ in \cite{BFM2014arX} we introduce with a new parameter $\alpha = \alpha(m)$ the new quantity $S(\alpha)$, where $\alpha$ will be chosen relatively small (we will see that $\alpha(m) = 1/(5m)$ is a good choice, for example).
Thus, let with a further parameter $\beta$
\begin{align}
\label{eq:3.7}
S(\alpha, \beta) &= \!\! \sum_{N < n \leqslant 2N}\! \biggl(\sum_{i = 1}^k 1_{\mathbb P} (n + h_i) - \beta m - \alpha \sum_{j = 1}^{4m + 1} \sum_{\substack{h, h'\in \mathcal H_j\\
h \neq h_2'}} 1_{\mathbb P}(n + h) 1_{\mathbb P}(n + h)'\biggr)\\
&\qquad \times \biggl(\sum_{\substack{d_1, \dots, d_k\\
d_i \mid n + h_i \forall i}} \lambda_{d_1, \dots, d_k}\biggr)^2
\nonumber
\end{align}
where under the summation sign we consider unordered pairs $h, h'\in \mathcal H_j$.
Let
\beq
\label{eq:3.8}
\beta = \beta(\alpha) = \max_{\ell \in \mathbb Z^+} \left(\ell - \alpha{\ell\choose 2}\right).
\eeq

Then the contribution of any set $\mathcal H_j$ to $S(\alpha, \beta) = S(\alpha)$ is at most $\beta$, so if we have for every $n \in (N, 2N]$ at most $m$ sets of the form $\mathcal H_j$ with
\beq
\label{eq:3.9}
\sum_{h \in \mathcal H_j} 1_{\mathbb P}(n + h) > 0,
\eeq
then consequently
\beq
\label{eq:3.10}
S(\alpha) \leqslant 0.
\eeq

In contrary to the choice $\varrho \in (0,1)$ and $\delta \varrho \log k = 2m$
of \cite{BFM2014arX} we will choose now $\delta \varrho \log k$ much larger
\beq
\label{eq:3.11}
\delta\varrho \log k = u: = \frac{4m + 1}{4\alpha}, \ \ \
\alpha = \frac{1}{5m}, \ \ \ \varrho\in(0,1).
\eeq
This implies with an easy calculation
\beq
\label{eq:3.12}
\beta = \frac{5m + 1}{2}.
\eeq
Using the same argument for the estimation of the negative double sum as \cite{BFM2014arX} we obtain a choice of a function $F$ such that
\begin{align}
\label{eq:3.13}
& S(\alpha) \!= \!\frac{N}{W} B^{-k} I_k(F) \biggl(\sum_{i = 1}^k \frac{u}{k} 
(1\! +\! O(\gamma))\! - \! \beta
m  - 4\alpha  \sum_{j = 1}^{4m + 1} \sum_{\substack{h, h'\in \mathcal H_j\\
h \neq h'}}\! \frac{u^2}{k^2} (1\! +\! O(\delta\! +\! \gamma))\!\biggr)\\
&= \! \frac{N}{W} B^{-k} I_k(F)\! \left(\!u(1\! +\! O(\gamma)) \! - \frac{(5 m\! +\! 1)m}{2} - 
\frac{4(4m\! + \! 1)}{5m} {k / (4m\! +\! 1) \choose 2}\! \frac{u^2}{k^2} (1 \! +\! O(\delta \! + \! \gamma))\!\right)\!.
\nonumber
\end{align}

By the above choice of the parameters in \eqref{eq:3.11} we have from \eqref{eq:3.13} with $\gamma = (\log k)^{-1/2}$
\begin{align}
\label{eq:3.14}
\frac{S(\alpha)WB^k}{NI_k(F)} &\geqslant  \frac{5m(4m + 1)(1 + O(\gamma))}{4} -
\frac{(5m + 1)m}{2}\\
&\qquad  - \frac{5m(4m + 1)(1 + O(\delta + \gamma))}{8}\nonumber \\
&= \frac{m(1 + O(m(\delta + \gamma)))}{8} > 0,
\nonumber
\end{align}
which contradicts to \eqref{eq:3.10}.

In order to see the validity of the last inequality we can choose
\beq
\label{eq:3.15}
m < (\log k)^{1/4} \Longleftrightarrow \delta \asymp \left(\frac{m^2}{\log k}\right)
\eeq
which implies $m\gamma = o(1)$ and $m\delta = o(1)$.
This contradiction proves our Theorem~\ref{th:1}.
Corollary~\ref{cor:2} follows from it in the same way as Corollary 1.2 from Theorem 1.1 in \cite{BFM2014arX}.

\noindent
J\'anos Pintz\\
R\'enyi Mathematical Institute\\
of the Hungarian Academy of Sciences\\
Budapest, Re\'altanoda u. 13--15\\
H-1053 Hungary\\
e-mail: pintz.janos@renyi.mta.hu

\end{document}